\documentclass{amsart}

\theoremstyle{plain}
\newtheorem{Prop}{Proposition}[section]
\newtheorem{Thm}[Prop]{Theorem}
\newtheorem{Cor}[Prop]{Corollary}

\newtheorem{Lem}[Prop]{Lemma}

\theoremstyle{definition}
\newtheorem{Def}[Prop]{Definition}

\theoremstyle{remark}
\newtheorem{Rem}[Prop]{Remark}
\newtheorem{Ex}[Prop]{Example}
\newtheorem{Problem}[Prop]{\bf Problem}

\def\dim{\mathop{\roman{dim}}}
\def\int{\mathop{\roman{int}}}

\def\1{^{-1}}

\def\dim{\text{dim}}
\def\asdim{\text{asdim}}

\def\T2{{\mathbf T_2}}
\def\NN{{\mathbb N}}

\def\UU{{\mathcal U}}
\def\LL{{\mathcal L}}
\def\VV{{\mathcal V}}

\def\RR{{\mathbb R}}
\def\ZZ{{\mathbb Z}}

\errorcontextlines=0 \numberwithin{equation}{section}

\begin{document}

\title[
Assouad-Nagata dimension of locally finite groups
and asymptotic cones]{Assouad-Nagata dimension of locally finite groups and asymptotic cones}
\author{J.~Higes}
\address{Departamento de Geometr\'{\i}a y Topolog\'{\i}a,
Facultad de CC.Matem\'aticas. Universidad Complutense de Madrid.
Madrid, 28040 Spain}
\email{josemhiges@yahoo.es}


\date{February 23, 2008}
\thanks{ The author is supported by Grant AP2004-2494 from the Ministerio de Educaci\' on y Ciencia, Spain and project MEC, MTM2006-0825. He also thanks Jerzy Dydak and J.M. Rodriguez Sanjurjo
for their support.}

\begin{abstract}In this paper we study two problems concerning Assouad-Nagata dimension:
\begin{enumerate}
\item Is there a metric space of positive Assouad-Nagata dimension such that all of its asymptotic cones are of Assouad-Nagata dimension zero? (Question 4.5 of \cite{Dydak-Higes})
\item Suppose $G$ is a locally finite group with a proper left invariant metric $d_G$. If $\dim_{AN}(G, d_G)>0$, is $\dim_{AN} (G, d_G)$ infinite?
(Problem 5.3 of \cite{Brod-Dydak-Lang})
\end{enumerate} 
The first question is answered positively. We provide examples 
of metric spaces of positive(even infinite) Assouad-Nagata dimension such that all of its asympotic cones are ultrametric. 
The metric spaces can be groups with proper left invariant metrics. \par
The second question has a negative solution. We show that for each $n$ there exists a locally finite group of Assouad-Nagata dimension $n$. 
As a consequence this solves for non finitely generated countable groups the question about 
the existence of metric spaces of finte asymptotic dimension whose asymptotic Assouad-Nagata dimension is larger but finite.

\end{abstract}
\maketitle
\tableofcontents
\section{Introduction}\label{Introduction}
Nagata dimension (also called Assouad-Nagata dimension)  was
introduced by Assouad in \cite{Assouad} influenced from the
papers of Nagata. This notion and its asymptotic version (asymptotic Assouad-Nagata dimension) has been studied
in recent years
as a geometric invariant, see for example
 \cite{Brod-Dydak-Higes-Mitra Lipschitz}, \cite{Dran-Smith2}, \cite{Dran-Zar} and specially \cite{Lang}.
In those papers many properties of the asymptotic dimension (see \cite{Bell-Dran} for a good survey about 
it) were generalized to asymptotic Assouad-Nagata dimension. 
\par
One interesting problem is about the relationship between the
Assouad-Nagata dimension of a metric space and the topological dimension of its asymptotic cones. \par
In \cite{Dydak-Higes} the following improvement of a result of \cite{Dran-Smith2} was obtained:
\begin{Thm}\label{DimensionOfCones} [Dydak, Higes \cite{Dydak-Higes}]
$dim(Cone(X, c, d) \le dim_{AN}(Cone(X, c, d) \le asdim_{AN}(X, d_X)$
for any metric space $(X, d_X)$.
\end{Thm}
The main idea of the proof of this theorem is that the $n$-dimensional Nagata property of a metric space is inherited
by its asymptotic cones. Recall that a metric space has the $0$-dimensional Nagata property if and only if it is ultrametric.
In section \ref{ConesAndQuestion} we show 
that the converse is not true by giving a class of non-ultrametric spaces such that all of its asymptotic cones are ultrametic. Such construction will also answer 
the first main question. \par 
Another problem that remains open is about the size
of the difference between the asymptotic dimension and the
asymptotic Assouad-Nagata dimension in a discrete group. P.Nowak
in \cite{Nowak} found for each $n \ge 1$ a finitely generated
group of asymptotic dimension $n$ but infinite Assouad-Nagata
dimension. If $n\ge 2$ such group can be finitely presented. As
a final problem of his paper, Nowak asked about the existence of
discrete groups such that the asymptotic dimension differs from
asymptotic Assouad-Nagata dimension but both of them are finite. The paper
of Nowak was complemented by one of Brodskiy, Dydak and Lang
\cite{Brod-Dydak-Lang} who related the growth of a finitely
generated group $G$ with the Assouad-Nagata dimension of the
wreath product $H\wr G$ with $H\ne 1$ a finite group.  As a
consequence of their results many examples of countable locally
finite groups with infinite Assouad-Nagata dimension can be found.
Countable locally finite groups satisfy many remarkable geometric
properties: they are the unique countable groups of
asymptotic dimension zero \cite{Smith}. Each metric space of
asymptotic dimension zero and bounded geometry can be embedded
coarsely in any infinite countable locally finite group. And they are the unique countable groups
that admits a proper left invariant ultrametric \cite{Brod-Dydak-Higes-Mitra}.
\par The main
target of section \ref{LocallyFiniteGroups} is to study asymptotic
Assouad-Nagata dimension of locally finite groups and countable
groups. In particular we find for each $n$ a locally finite group
$(G, d_G)$ with $d_G$ a proper left invariant metric that
satisfies $asdim_{AN} (G, d_G) = n$. In \cite{Brod-Dydak-Levin-Mitra} a method was shown to build a locally finite group 
of infinite Assouad-Nagata dimension. Our construction seems quite similar to that one but it was not clear how to define the group and the metric in such a way its asymptotic Assouad-Nagata dimension 
was finite and positive. This problem was asked explictly by two of the authors of \cite{Brod-Dydak-Levin-Mitra} and Lang in \cite{Brod-Dydak-Lang}. 
We also show in Corollary \ref{GeneralizationToCountable} that for each
$n\ge 0$ and $k\ge 0$ there is a countable group $G$ and a proper
left invariant metric $d_G$ such that $G$ is of asymptotic
dimension $n$ but $asdim_{AN}(G, d_G) = n +k$. This solves Nowak's
problem for countable groups and answers negatively the second main problem of this paper.
As far as we know the problem for finitely generated groups and finitely presented groups remains open.

\section{Asymptotic cones and cubes}\label{AsymptoticConesAndAssouadNagataDimension}

Let $s$ be a positive real number. An {\it $s$-scale chain} (or $s$-path)
 between two points $x$ and $y$ of a metric space $(X, d_X)$ is defined as
a finite sequence points
$\{x= x_0, x_1, ..., x_m = y\}$ such that
$d_X(x_i, x_{i+1}) < s$ for every $i = 0, ..., m-1$.
A subset $S$ of a metric space $(X, d_X)$ is said to be {\it $s$-scale connected} if there exists an $s$-scale chain
contained in $S$ for every two elements of $S$.
\begin{Def} A metric space $(X, d_X)$ is said to be of
{\it asymptotic dimension} at most $n$ (notation $\asdim(X, d) \le
n$) if there is an increasing function $D_X: \RR_+ \to \RR_+$ such
that for all $s> 0$ there is a cover $\UU =\{\UU_0, ...,\UU_n\}$
so that the $s$-scale connected components of each $\UU_i$ are
$D_X(s)$-bounded i.e. the diameter of such components is bounded
by $D_X(s)$. \par The function $D_X$ is called an {\it
$n$-dimensional control function} for $X$. Depending on the type
of $D_X$ one can define the following two invariants: \par A
metric space $(X, d_X)$ is said to be of {\it  Assouad-Nagata
dimension} at most $n$ (notation $\dim_{AN} (X, d) \le n$) if it
has an $n$-dimensional control function $D_X$ of the form
 $D_X(s) = C\cdot s$ with $C>0$ some fixed constant. \par
A metric space $(X, d_X)$ is said to be of
{\it asymptotic Assouad-Nagata dimension} at most $n$ (notation $\asdim_{AN} (X, d) \le n$) if it has an
$n$-dimensional control function $D_X$ of the form
 $D_X(s) = C\cdot s +k$ with $C >0$ and $k\in \RR$ two fixed constants.
\end{Def}
One important fact about the asymptotic dimension is that it is invariant under coarsely equivalences. Given a map $f: (X, d_X) \to (Y, d_Y)$ between two metrics spaces
it is said to be a {\it coarse embedding} if there exists two increasing functions $\rho_+: \RR_+ \to \RR_+$ and $\rho_-: \RR_+\to \RR_+$ with 
$\lim\limits_{x\to \infty} \rho_-(x) =\infty$
such that:
\[\rho_-(d_X(x, y)) \le d_Y(f(x), f(y))\le \rho_+(d_X(x, y)) \text{ for every } x, y \in X.\]
The functions $\rho_+$ and $\rho_-$ are usually called {\it contraction and dilatation} functions of $f$ respectively. \par
Now a {\it coarse equivalence} between two metrics spaces $(X, d_X)$ and $(Y, d_Y)$ is defined as a coarse embedding $f: (X, d_X)\to (Y, d_Y)$ for which there exists a constant 
$K >0$ such that
$d_Y(y, f(X))\le K$ for every $y \in Y$. If there exists a coarse equivalence between $X$ and $Y$ both spaces are said to be {\it coarsely equivalent}. \par
Next result relates $n$-dimensional control functions and coarse embeddings. 

\begin{Prop}\label{Traductor} Let $(X, d_X)$ and $(Y, d_Y)$ be two metric spaces and let $f: (X, d_X) \to (Y, d_Y)$ be a coarse embedding with $\rho_+$ and $\rho_-$ the contraction 
and dilatation functions of $f$ respectively. If $D_Y^n$ is an $n$-dimensional control function of $(Y, d_Y)$ then the funcion $D_X^n$ defined by 
$D_X^n = \rho_-^{-1} \circ D_Y^n \circ \rho_+$ is an $n$-dimensional control function of $(X, d_X)$. 
\end{Prop}
\begin{proof}
Fix $s >0$ a positive real number. As $D_Y^n$ is an $n$-dimensional control function there exists a cover $\UU =\{\UU_0, ...,\UU_n\}$ in $Y$ so that the $\rho_+(s)$-scale connected components 
of each $\UU_i$ are $D_Y^n(\rho_+(s))$-bounded. Take the cover $\VV = \{\VV_0, ..., \VV_n\}$ in $X$ defined as $\VV_i = f^{-1}(\UU_i)$. Notice that if two points $x, y \in X$ 
satisfies $d_X(x, y) < s$ then $d_Y(f(x), f(y)) < \rho_+(s)$. Hence given an $s$-scale chain $\{x_0, x_1,..., x_m\}$ in $X$ we get that $\{f(x_0), f(x_1), ..., f(x_m)\}$ is an 
$\rho_+(s)$-scale chain. Therefore $d_Y(f(x_0), f(x_m)) \le D_Y^n(\rho_+(s))$ what implies $d(x_0, x_m) \le \rho_-^{-1}(D_Y^n(\rho_+(s)))$ and 
$D_X^n = \rho_-^{-1} \circ D_Y^n \circ \rho_+$ is an $n$-dimensional control function.   
\end{proof}
The following easy corollary will be useful in the next section:

\begin{Cor} \label{PolynomialCorollary} 
If $(X, d)$ is a metric space and $asdim_{AN}(X, log(1+d))\le n$ then there is a polynomial $n$-dimensional control function of $(X, d)$.
\end{Cor}
\begin{proof} 
Suppose $\asdim_{AN}(X, log(1+d)) \le n$.  This implies 
there exists a linear $n$-dimensional control function $D_X^{n}(s) = C\cdot s + b$ of $(X, log(1+d)$ with $C$ and $b$ two fixed constants. 
Suppose without lost of generality that $C\in \NN$. 
Now the identity $(X, d|_x) \to (X, log(1+d))$ is clearly a coarse equivalence with $\rho_+(d) = log(1+d) = 
\rho_-(d)$. By proposition \ref{Traductor} we will get that $Q_X^{n} = \rho_-^{-1}\circ D_X^n \circ \rho_+ $ is an 
$n$-dimensional control function of $(X, d)$. That means $Q_X^n(s) = (10^b\cdot (1+s)^C - 1$. Therefore there is a polynomial dimensional control function of $(X,d)$.        
\end{proof}
It is clear that the asymptotic dimension of a metric space is
less than or equal the asymptotic Assouad-Nagata dimension and this is
greater or equal than the Assouad-Nagata dimension. In
\cite{Brod-Dydak-Higes-Mitra Lipschitz} it was shown that for a
discrete space the asymptotic Assouad-Nagata dimension and the
Assouad-Nagata dimension are equal.\par 

Our target in this
section is to find some sufficient conditions that give lower
bounds for the Assouad-Nagata dimension. Next definition plays an
important role for such purpose.

\begin{Def}
We define an {\it $n$-dimensional dilated cube}
 in a metric space $(X, d_X)$ as a dilatation function $f:\{0, 1, ..., k\}^n\to X$,
that means there exists a constant $C \ge 1$
(dilatation constant) such that
$C\cdot \|x - y\|_1 = d_X(f(x), f(y))$ for every $x, y \in \{0, 1,...,k\}^n$.
\end{Def}
\begin{Rem}
$n$-dimensional dilated cubes are particular cases of the
$n$-dimensional $s$-cubes introduced by Brodskiy, Dydak and Lang
in \cite{Brod-Dydak-Lang}. Recall that an $n$-dimensional $s$-cube
in a metric space $(X, d_X)$ is defined as a function $f:\{0,
1,..., k\}^n\to X$ with the property that $d(f(x),f(x+e_i))< s$
for all $x \in \{0, 1, ..., k\}^n$ such that $x+e_i\in \{0,
1,...,k\}^n$ with $e_i$ belonging to the standard basis of
$\RR^n$.\par
In the whole paper we will take in $\RR^n$ the $l_1$-metric instead of the euclidean metric.
\end{Rem}

Now we will relate the existence of some sequences of
$n$-dimensional dilated cubes in a space with  the existence
of cubes of the form $[0, s]^n \subset \RR^n$ in its asymptotic cones. Let $(X, d_X)$ be a
metric space. Given a non-principal ultrafilter $\omega$ of $\NN$
and a sequence $\{x_n\}_{n\in\NN}$ of points of $X$, the {\it
$\omega$-limit} of $\{x_n\}_{n\in\NN}$ (notation:
$\lim\limits_{\omega} x_n$)is defined as the element $y$ of $X$
such that for every neighborhood $U$ of $y$ the set $F_{U} = \{n|
x_n \in U\}$ belongs to $\omega$. Analogously if for every ball
$B(x, r)$ of radius $r$ the set $F_{B(x,r)} = \{n | x_n \in X
\setminus B(x, r)\}$ belongs to $\omega$ then it is said that the
$\omega$-limit of $\{x_n\}_{n\in \NN}$ is infinity and the
sequence is an $\omega$-divergent sequence. It can be proved
easily that the $\omega$-limit always exists in a compact space.
\par Assume $\omega$ is a non principal ultrafilter of $\NN$.  Let
$d =\{d_n\}_{n\in\NN}$ be an
$\omega$-divergent sequence of positive real numbers  and 
let $c = \{c_n\}_{n \in \NN}$ be  any
sequence of elements of $X$. Now we can construct the
{\it asymptotic cone} (notation: $Cone_{\omega}(X, c, d)$) of $X$ as follows:
\par Firstly define the set of all sequences $\{x_n\}_{n\in\NN}$
of elements of $X$ such that $\lim\limits_{\omega}\frac{d_X(x_n,
c_n)}{d_n}$ is bounded. In such set take the pseudo metric given
by: \[D(\{x_n\}_{n\in\NN}, \{y_n\}_{n\in\NN}) =
\lim\limits_{\omega}\frac{d_X(x_n, y_n)}{d_n}.\] By identifying
sequences whose distances is $0$ we get the metric space
$Cone_{\omega}(X, c, d)$.\par Asymptotic cones were originally
introduced by Gromov in \cite{Gromov}. There has been a lot of
research relating properties of groups with topological properties
of its asymptotic cones. For example a finitely generated group is
virtually nilpotent if and only if all its asymptotic cones are
locally compact \cite{Gromov2} or a group is hyperbolic if and only if all of its asymptotic cones are $\RR$-trees (\cite{Gromov} and \cite{Drutu}). 

\begin{Prop}\label{SequenceOfCubesImpliesCube} Let $(X, d_X)$ be a
metric space and let $\{f_m\}_{m =1}^{\infty}$ be a sequence of
$n$-dimensional dilated cubes, $f_m: \{0, 1,..., k_m\}^n \to X$
such that $\lim\limits_{\omega} k_m = \infty$ for some non principal ultrafilter $\omega$ of $\NN$.
If $\{d_m\}_{m =1}^{\infty}$ is an $\omega$-divergent sequence of positive real numbers that satisfies
\[\lim\limits_{\omega}
\frac{C_m \cdot k_m}{d_m} = s\] with $\{C_m\}_{m= 1}^{\infty}$
the sequence of dilatation constants and $ 0 \le s < \infty$,  then $[0, s]^n \subset Cone_{\omega}(X, c, d)$ where
$c = \{f_m(0)\}_{m =1}^{\infty}$.
\end{Prop}
\begin{proof}

Let us prove firstly the case $n = 1$.  For each $t \in [0, s]$ let $A_m$ be the subset of elements of 
$\{0, 1,..., k_m\}$ such that, for every $x \in A_m$, $\frac{C_m\cdot x}{d_m}$ is closest to $t$. It means $C_m\cdot x$ is closest to $d_m \cdot t$. Take now
the sequence $\{r_m^t\}_{m = 1}^{\infty}$  where $r_m^t$ is the
minimum element of $A_m$.\par
Define the map $g:[0, s] \to Cone_{\omega}(X, c, d)$ by $g(t) = x^t$
if the sequence $\{f_m(r_m^t)\}_{m = 1}^{\infty}$ is in the class $x^t$. As:
\[\lim\limits_{\omega} \frac{d(f_m(0), f_m(r_m^t))}{d_m} =
\lim\limits_{\omega} \frac{C_m\cdot r_m^t}{d_m} \le \lim\limits_{\omega} \frac{C_m \cdot k_m}{d_m} = s\]
the map is well defined. Let us prove that it is in fact an
isometry. From the definition of $r_m^t$ we get that if $t_1 < t_2$ then
$r_m^{t_1}\le r_m^{t_2}$ what implies $\lim\limits_{\omega}\frac{d(f_m(r_m^{t_1}), f_m(r_m^{t_2}))}{d_m}
= \lim\limits_{\omega} \frac{C_m(r_m^{t_2}-r_m^{t_1})}{d_m}$.
So the unique thing we need to show is that $\lim\limits_{\omega} \frac{C_m \cdot r_m^t}{d_m} = t$
for every $t$. Firstly notice that as $\lim\limits_{\omega} \frac{C_m \cdot k_m}{d_m} = s$ and $s \ge t$ then there exists 
an $F \in \omega$ such that $\frac{C_m \cdot k_m}{d_m} \ge t$ for every $m \in F$. We have also $\lim\limits_{\omega}\frac{C_m}{d_m} = 0$ as 
$\lim\limits_{\omega}k_m = \infty$ but $\lim\limits_{\omega} \frac{C_m \cdot k_m}{d_m} = s < \infty$. This implies that  
given $\epsilon >0$ there exists $G_{\epsilon} \in \omega$ such that $\frac{C_m}{d_m}<\epsilon$ for every $m \in G_{\epsilon}$. Therefore 
if $m \in F \cap G_{\epsilon}$ we have $|C_m\cdot r_m^t -d_m\cdot t|<  C_m$ and then $|\frac{C_m \cdot r_m^t}{d_m} - t|< \frac{C_m}{d_m}\le \epsilon$.

Now let us do the general case. Let $(s_1, ..., s_n) \in [0, s]^n$. By the previous case we get that for every $j = 1,..., n$
 there exists a sequence $\{r_m^{s_j}\}_{m \in \NN}$ with $r_m^{s_j}\in \{0,1,...,k_m\}$ such that 
$\lim\limits_{\omega} \frac{C_m\cdot r_m^{s_j}}{d_m} = s_j$. In a similar way as before we construct a map $g: [0,s]^n \to Cone_{\omega}(X, c, d)$ 
by defining $g(s_1, ..., s_m)$ as the class that contains the sequence $\{f_m(r_m^{s_1}, ..., r_m^{s_n})\}_{m = 1}^{\infty}$. 
To finish the proof it will be enough to check that
 for every $s, t \in [0, s]^n$ with $s = (s_1,..., s_n)$ and $t = (t_1,..., t_n)$, the following equality holds:
\[   \lim\limits_{\omega}\frac{d_X(f_m(r_m^{s_1},..., r_m^{s_n}), f_m(r_m^{t_1},..., r_m^{t_n}))}{d_m} = \sum_{i= 1}^n |s_i -t_i|\]
As $f_m$ is a dilatation of constant $C_m$ we can write:
\[\lim\limits_{\omega}\frac{d_X(f_m(r_m^{s_1},..., r_m^{s_n}), f_m(r_m^{t_1},..., r_m^{t_n}))}{d_m} = \sum_{i= 1}^n \lim\limits_{\omega} 
\frac{C_m\cdot |r_m^{s_i}- r_m^{t_i}|}{d_m}\]
And again by the case $n =1$ we can deduce that the last term satisfies the equality: 
 \[\sum_{i= 1}^n \lim\limits_{\omega} 
\frac{C_m\cdot |r_m^{s_i}- r_m^{t_i}|}{d_m} = \sum_{i=1}^n |s_i-t_i|\].

\end{proof}

Combining theorem \ref{DimensionOfCones} with proposition \ref{SequenceOfCubesImpliesCube}
we can get the following lower bound of Assouad-Nagata dimension for certain spaces.
Such estimation will be very useful to prove the main results of section \ref{LocallyFiniteGroups}.

\begin{Cor}\label{SequenceImpliesNagata}
If a metric space $(X, d_X)$ contains a sequence $\{f_m\}_{m=1}^{\infty}$ of
$n$-dimensional dilated cubes $f_m:\{0, 1,..., k_m\}^n \to X$ with
$\lim\limits_{m\to \infty} k_m = \infty$ then $asdim_{AN}(X, d_X) \ge n$.
\end{Cor}
\begin{proof}
Let $\{C_m\}_{m\in \NN}$ be the sequence of dilatation constants.
Define the divergent sequence $d = \{d_m\}_{m \in \NN}$ as $d_m = C_m\cdot k_m$.
Then for any non principal ultrafilter $\omega$ of $\NN$ the hypothesis of proposition
\ref{SequenceOfCubesImpliesCube} are satisfied with $s = 1$ so we get $[0, 1]^n \subset Cone_{\omega}(X, c, d)$.
Applying \ref{DimensionOfCones} we obtain immediately:
\[n \le dim(Cone_{\omega}(X, c, d))\le dim_{AN}(Cone_{\omega}(X, c, d)) \le asdim_{AN}(X, d_X)\]
\end{proof}

\section{Ultrametric Asymptotic cones} \label{ConesAndQuestion}
The main aim of this section is to find metric spaces of positive asymptotic Assouad Nagata dimension such that all of its asymptotic cones are  ultrametric spaces. 
In particular we are interested when the space is a group with a proper left invariant metric. 
A metric $d_G$ defined in a group $G$ is said to be a
{\it proper left invariant metric} if it satisfies the following conditions:
\begin{enumerate}
\item $d_G(g_1\cdot g_2, g_1\cdot g_3) = d_G(g_2, g_3)$ for every $g_1, g_2, g_3 \in G$.
\item For every $K>0$ the number of elements $g$ of $G$ such that $d(1, g)_G \le K$ is finite.
\end{enumerate}
We will say that a metric $d_X$ defined in a set $X$ is an {\it asymptotic ultrametric}
if there exists a constant $k \ge 0$ such that:
\[d_X'(x, y) \le \max\{d_X'(x, z), d_X'(y,z)\}+k \text{ for every } x, y, z \in X\]
If $k = 0$ the space is said to be {\it ultrametric}. \par
Next theorem is a combination of two different results of \cite{Brod-Dydak-Higes-Mitra Lipschitz} and \cite{Brod-Dydak-Levin-Mitra}.

\begin{Thm}\label{AsymptoticUltrametric} For every metric space $(X,d_X)$ of finite asymptotic dimension there exists a metric $d_X'$ coarsely equivalent to $d_X$ such that:
\begin{enumerate}
\item $asdim_{AN} (X, d_X') = asdim(X, d_X)$.
\item $d_X'$ is an asymptotic ultrametric.
\end{enumerate}
Moreover if $(X, d_X)$ is a countable group with $d_X$ a proper left invariant metric we can take $d_X'$ as a proper left invariant metric.
\end{Thm}

\begin{proof}
The first part of the theorem follows directly from the proof of  [\cite{Brod-Dydak-Higes-Mitra Lipschitz}, Theorem 5.1.]\par
The second part can be got by an easy modification of the proof of [\cite{Brod-Dydak-Levin-Mitra}, proposition 6.6. 
\end{proof}
Next result appeared in Gromov [\cite{Gromov}, page] but without proof. We provide a proof here.
\begin{Prop}\label{EpsilonKUltrametric}Let  $(X, d_X)$ be a metric space that satisfies:
\[d_X(x, y) \le (1+ \epsilon) \cdot \max\{d_X(x, z), d_X(y, z)\} +k\text{ for every } x, y , z \in X\]
where $\epsilon$ is some function in $d = d_X(x, y)$ which goes to 0 when $d$ goes to $\infty$ and $k>0$ some fixed constant.   
Then every asymptotic cone of $(X, d_X)$ is an ultrametric space
\end{Prop}
\begin{proof}
Suppose $(X, d_X)$ is an asymptotic ultrametric space with constant $k>0$.
Let $x$, $y$ and $z$ be three points of  $Cone_{\omega}(X, c, d)$ and
let  $\{x_n\}_{n=1}^{\infty}, \{y_n\}_{n=1}^{\infty}, \{z_n\}_{n=1}^{\infty}\subset X$ be three sequences in the classes $x$, $y$ and $z$ respectively. Without loss of generality 
assume $D(x,y) \ge D(x, z)$ and $D(x, y) \ge D(y,z)$. It will be enough to check that $D(x, y) \le \max\{D(x, z), D(y, z)\}$. 

First notice that $\lim_{\omega}\epsilon(d_X(x_n, y_n)) = 0$. If not there exists an $M>0$ such that the set $H =\{n | \epsilon(d_X(x_n, y_n)) \ge M\}$ is in $\omega$. But 
if we assume the non trivial case $D(x, y) \ne 0$ then the set $G = \{d_X(x_n, y_n)| n \in H\}$ is unbounded. Taking now a divergent subsequence 
$\{d_X(x_{n_i}, y_{n_i})\}_{i = 1}^{\infty}\subset G$ we get a contradiction as $\lim_{i \to \infty} \epsilon(d_X(x_{n_i}, y_{n_i}) = 0$ but $d_X(x_{n_i}, y_{n_i})\ge M $ for 
every $i \in \NN$. 

Take now $F$ the subset of natural numbers defined as:
\[F =\{n | d_X(x_n, z_n) \ge d_X(y_n, z_n)\}\] As $\omega$ is an ultrafilter then $F\subset \omega$ or $\NN \setminus F \subset \omega$.
Let us assume the first case. In this first case  we will obtain that $D(x, z) \ge D(y, z)$ and by
the asymptotic ultrametric property we deduce that:
\[\frac{d_X(x_n, y_n)}{d_n}\le \frac{(1+\epsilon(d_X(x_n, y_n)))\cdot d_X(x_n, z_n)}{d_n}+ \frac{k}{d_n}\text{ for every } n\in F. \]
Taking limits in the previous inequality and applying the fact $\lim_{\omega}\epsilon(d_X(x_n, y_n)) = 0$ 
we obtain $D(x, y)\le D(x,z) = \max\{D(x,z), D(y, z)\}$. 
Finally let us do 
the case $\NN \setminus F \subset \omega$. This implies 
the set $\NN \setminus F  = \{n | d_X(y_n, z_n) > d_X(x_n, z_n)\}$ belongs to the ultrafilter what yields $D(y, z) \ge D(x, z)$ appying now the same reasoning we get:
 \[\frac{d_X(x_n, y_n)}{d_n}\le \frac{(1+ \epsilon(d_X(x_n, y_n)))\cdot d_X(y_n, z_n)}{d_n}+ \frac{k}{d_n}\text{ for every } n\in \NN \setminus F. \]
Taking again limits we obtain $D(x, y)\le D(y,z) = \max\{D(x,z), D(y, z)\}$ what finish the proof. 
\end{proof}

\begin{Cor}\label{ConesOfAsymptoticUltrametricAreUltrametric} If $(X,d_X)$ is a metric space
with $d_X$ an asymptotic ultrametric then $Cone_{\omega}(X,c, d)$ is an ultrametric space for every $c$, $d$ and $\omega$.
\end{Cor}

Recall that ultrametric spaces are of Assouad-Nagata dimension zero. So 
an immediate consequence of  theorem \ref{AsymptoticUltrametric}
and proposition \ref{ConesOfAsymptoticUltrametricAreUltrametric}
is the following theorem that solves the first target of this paper:

\begin{Thm}For every metric space $(X, d_X)$ with finite asymptotic dimension there exists a coarsely equivalent
metric $d_X'$ such that $asdim_{AN}(X, d_X') = asdim(X, d_X)$ and
every asymptotic cone of $(X, d_X')$ is an ultrametric space.\par
Moreover if $(X, d_X)$ is a countable group with $d_X$ a proper left invariant metric then we can take $d_X'$ a
proper left invariant metric.
\end{Thm}

We could ask now about the existence of a metric space of infinite asymptotic Assouad-Nagata dimension such that all its asymptotic 
cones are of Assouad-Nagata dimension zero, in particular ultrametric spaces. The rest of the section is devoted to show an example of this type. 
Such example will be a finitely generated group with some proper left invariant metric. We will need the following result of \cite{Brod-Dydak-Lang}.       

\begin{Thm}(Brodskiy, Dydak and Lang \cite{Brod-Dydak-Lang})\label{BrodDydakLang} Suppose $H$ and $G$ are finitely generated and $K$ is the kernel of the projection of $H\wr G\to G$ equiped with the metric 
induced from $H\wr G$. If $\gamma $ is the growth function of $G$ and $D_K^{n-1}$ is an $n-1$ dimensional control function of $K$, then the integer part of $\frac{\gamma(r)}{n}$
is at most $D_K^{n-1}(3\cdot n \cdot r)$.
\end{Thm}

\begin{Thm}\label{InfiniteNagataUltrametric} Let $G$ be a finitely generated group of exponential growth and let $H$ be a finite group. Suppose $d$ is a word metric of $H \wr G$ then:
\begin{enumerate}
\item All the asymptotic cones of $(H\wr G, log(1+d))$ are ultrametric. 
\item $\asdim_{AN} (H\wr G, log(1+d)) = \infty$.
\end{enumerate}
\end{Thm}
\begin{proof}
The first assertion is a consequence of propositon \ref{EpsilonKUltrametric}. In fact it is not hard to check that if the metric  $D(x, y)$ of $X$ is of the form $D(x, y) = log(1+d(x,y))$ where $d(x, y)$ is another distance, 
then $\epsilon(d) = \frac{log(2)}{log(\frac{d}{2}+1)}$ satisfies the conditions of the cited proposition.\par
Second assertion is a consequence of corollary \ref{PolynomialCorollary} and theorem \ref{BrodDydakLang}. 
The proof is by contradiction. Suppose $\asdim_{AN}(H\wr G, log(1+d)) \le n$. By corollary \ref{PolynomialCorollary} we get that there is a polynomial 
$n$-dimensional control function 
of $(H\wr G, d)$ and then there is a polynomial $n$-dimensional control function of the kernel $K \subset H \wr G$ with the restricted metric. But as the growth of $G$ is exponential by 
theorem \ref{BrodDydakLang} any $n$-dimensional control function of $(K, d|_k)$ must be at least exponential, a contradiction.       
\end{proof}
\begin{Rem}
As the metric spaces $(H\wr G, d)$ and $(H\wr G, log(1+d))$ are coarsely equivalent then both have the same asymptotic dimension. Hence we can also assume that the group 
$(H \wr G, log(1+d))$ is of finite 
asymptotic dimension. Let us show an example of this fact.
\end{Rem}   
\begin{Ex}
Let 
$G = F_2$ be the free group of two generators and  $H = \ZZ_2$. For any word metric $d$ of $\ZZ_2 \wr F_2$ we have $\asdim (\ZZ_2 \wr F_2, d) = 1$ and 
such group satisfies the conditions of 
the theorem. 
\end{Ex}

It is clear that the kernel $K$ of the projection $H \wr G \to G$ is a locally finite group when $H$ is finite. From the proof of \ref{InfiniteNagataUltrametric} 
we get easily:

\begin{Cor} There exists a locally finite group $K$ with a proper left invariant metric $d_K$ such that $\asdim_{AN}(K, d_K) = \infty$  and all 
of its asymptotic cones are ultrametric.
\end{Cor}
 
\begin{Rem} Notice that if $G$ is a finitely generated group and $d_G$ is any word metric then every asympotic cone $Cone_{\omega}(G, c, d)$ is geodesic 
what implies  
$dim_{AN}(Cone_{\omega}(G, c, d) \ge 1$. Therefore the metrics $d_G'$ of the theorems of this section can not be quasi-isometric to any word metric of $G$.
\end{Rem}

\section{Locally finite groups and positive Assouad Nagata dimension}\label{LocallyFiniteGroups}
In this section we will study  the Assouad-Nagata dimension of locally finite groups.
Firstly let us check that all locally finite
groups admits a proper left invariant metric with positive Assouad-Nagata dimension.
\par

Given a countable group $G$ we want to build a proper left invariant metric on it. 
Associated to each proper left invariant metric there exists a proper norm. 
A map $\|\cdot\|_G: G \to \RR_+$ is called to be a {\it proper norm} if  it satisfies
the following conditions:
\begin{enumerate}
\item $\|g\|_{G} = 0$ if and only if $g$ is the neutral element of $G$.
\item $\|g\|_G = \|g^{-1}\|_G$ for every $g \in G$.
\item $\|g\cdot h\|_G \le \|g\|_G +\|h\|_G$ for every $g, h \in G$.
\item For every $K>0$ the number of elements of $G$ such that $\|g\|_G\le K$ is finite.
\end{enumerate}
So if we build a proper norm $\|\cdot\|_G$ in $G$ then the map $d_G(g, h) = \|g^{-1}\cdot h\|_G$ defines a proper left invariant metric. 
Conversely for every proper left invariant metric $d_G$ the map $\|g\|_G = d_G(1, g)$ defines a proper norm.\par
One method of obtaining a norm in a countable group was described
by Smith in \cite{Smith}(see also \cite{Shalom}). Let $S$ be a symmetric system of
generators(possibly infinite) of a countable group $G$ and let $w:
L \to \RR_+$ be a function({\it weight function}) that satisfies:
\begin{enumerate}
\item $w(s) = 0$ if and only if $s = 1_G$
\item $w(s) = w(s^{-1})$.
\end{enumerate}
Then the function $\|\cdot\|_w: G\to \RR_+$ defined by:
\[\|g\|_w = \min\{\sum_{i=1}^n w(s_i)| x = \Pi_{i=1}^n s_i, \text{ }s_i \in S\}\]
is a norm. Moreover if $w$ satisfies also that $w^{-1}[0, N]$ is finite for every $N$ then $\|\cdot\|_w$ is a proper norm.
\par Notice that if we define $w(g) = 1$ for
all the elements $g \in S$ of a finite generating system $S\subset G$ ($G$ a finitely generated group) we will obtain the usual word metric.
\par This method for finite groups
has the following nice(and obvious) extension property:
\begin{Lem} \label{ExtensionNorm} Let $G$ be a countable group and let $(G_1, d_{G_1})$ be a finite subgroup of $G$ with $d_{G_1}$ a proper left invariant metric. Let $S\subset G$ be a symmetric subset
of $G$ such that:
\begin{enumerate}
\item $S \cap G_1 = \emptyset$.
\item $G$ is generated by $G_1 \cup S$.
\end{enumerate}
If $\|\cdot\|_w: G \to \RR_+$ is a norm defined by a weight function $w: G_1 \cup S \cup S^{-1} \to \RR_+$ that satisfies:
\begin{enumerate}
\item $w(g) = \|g\|_{G_1}$ if $g \in G_1$
\item $w(g) \ge diam( G_1)$ if $g \in S \cup S^{-1}$.
\end{enumerate}
then for every $g\in G_1$ $\|g\|_w = \|g\|_{G_1}$.
 \end{Lem}
A group is said to be {\it locally finite} if all of its finitely
generated subgroups are finite.\par In a locally finite group $G$
we can take a filtration $\LL$ of finite subgroups $\LL = \{\{1\}
= G_0 \subset G_1 \subset G_2 ...\}$ of $G$.  Lemma
\ref{ExtensionNorm} applied successively to $\LL$ allows us to
build a sequence of norms $\|\cdot \|_i: G_i \to \RR_+$ and a norm
$\|\cdot\|_G: G\to \RR_+$ such that the restriction of
$\|\cdot\|_G$ to $G_i$ coincides with $\|\cdot\|_{G_i}$ and each
$\|\cdot\|_i$ is an extension of $\|\cdot\|_{i-1}$. This idea was
used in \cite{Brod-Dydak-Higes-Mitra} to prove that each locally
finite group is coarsely equivalent to a direct sum of cyclic
groups. \par Let $\LL = \{\{0\} = G_0\subset G_1\subset...\}$ be a
filtration of a countable group $G$. We will say that the sequence
$\{g_n\}_{n =1}^{\infty}$ of elements of $G$ is a {\it system of
generators} of $\LL$ if for every $i \ge 1$ there exists an $n_i$
such that $G_i = < g_1,..., g_{n_i}>$. If it is also satisfied
that $n_i = i$ and $G_{i-1}\ne G_i$ for every $i\ge 1$ then we
will say that the filtration $\LL$ is a {\it one-step ascending
chain}. In one-step ascending chains we can estimate easily the
cardinality of $G_i$.
 \begin{Lem}\label{IncreaseImpliesPower} Let $G$ be a locally finite group and let $\LL =\{\{1\} = G_0 < G_1 < ...\}$ be a one-step ascending chain then $|G_i| \ge 2^i$ for every $i \in \NN$.
\end{Lem}
\begin{proof}
Let $\{g_i\}_{i \in \NN}$ be a sequence of generators of $\LL$.
For $i = 1$ the result is obvious. Let us assume the result is true for some $i\ge 1$. Take $2^i$ different elements
$\{x_n\}_{n =1}^{2^i}$ of $G_i$. The subset $X$ of $G_{i+1}$ defined as
 $X = \{x_n\}_{n =1}^{2^i} \cup \{x_n\cdot g_{i+1}\}_{n = 1}^{2^i}$ contains $2^{i+1}$ different elements as $g_{i+1} \not\in G_i$.
\end{proof}

\begin{Thm} Let $G$ be a locally finite group. For every increasing function $f: \RR_+ \to \RR_+$ such that  $\lim\limits_{x \to \infty} f(x) = \infty$
there exists a proper left invariant metric $d_G$
in $G$ such that $f$ is not a $0$-dimensional control function.
\end{Thm}
\begin{proof}
The proof is by contradiction. Let $\LL = \{ G_0 < G_1<...\}$ be a
one-step ascending chain of $G$ with $\{g_i\}_{i=1}^{\infty}$ a sequence of
generators of $\LL$. Assume there exists some function $f$ as in the
hypothesis that is a $0$-dimensional control function such that for every proper left invariant metric $d_G$. In such case we
claim:\par {\it There exists a $K>0$ such that for every $i>0$
every element $g \in G_i$ has length less than or equal $K$ as a word
of $\{g_1, ..., g_i\}$.} \par
Suppose that the claim were false. In such case we will show there is a
proper left invariant metric $d_G$ such that $f$ is not a $0$-dimensional control function.
Fix
$J_1 = 1$ and take $h_1\in
G_{i_1}$ an element such that its minimum lenght in  $G_{i_1}$ is greater than
$f(1)$. Define the norms in the system of generators
$\{g_j\}_{j=1}^{i_1}$ as $\|g_j^{\epsilon}\|_{G_{i_1}} = 1$ with $\epsilon = \pm 1$. By viewing such norms as weights we have
a proper left invariant metric in $G_{i_1}$. Notice that
$\|h_1\|_{G_{i_1}} > f(J_1)$ and $h_1$ is in the same $1$-connected
component as the origin. Now suppose defined a proper left
invariant metric in $G_{i_r}$  that verifies the following property: \par
There exists an $h_r \in
G_{i_r}$ with $\|h_r\|_{G_{i_r}}
> f(J_{r})$, $J_{r} \ge diam(G_{i_{r-1}})$ and $h_r$ is in
the same $J_r$-connected component as the origin.
\par
Let $J_{r+1}> diam( G_{i_r})$ and let $h_{r+1}$ be
an element of the group $G_{i_{r+1}}$ with its minimum length greater than $f(J_{r+1})$.
Apply \ref{ExtensionNorm} to $G_{i_r} \cup S$ with $S = \{g_{i_r+1}, ...,g_{i_{r+1}}\}$ and a weight function $w$ that satisfies $w(g) = J_{r+1}$
for every generator $g  \in S \cup S^{-1}$. Then we get a norm in $G_{i_{r+1}}$ that is an extension of $\|\cdot\|_{G_{i_r}}$ with the same property as above. Repeating this
procedure we can get a proper norm $\|\cdot\|_G$
defined in $G$ so that for every $J_r$ there is an element $h_r$ with
$\|h_r\|_G > f(J_r)$ and $h_r$ is in the same $J_r$ connected
component as the origin. We deduce that $f$ can not be a $0$-dimensional control function of $d_G$ so the claim must be true. \par
Now applying an easy combinatorial argument  we can estimate the
cardinality of every subgroup $G_i$:
\[|G_i| \le \sum_{j= 0}^K(2\cdot i)^j\le (K+1)\cdot 2^K\cdot i^{K}\] This
contradicts lemma \ref{IncreaseImpliesPower} for $i$ sufficiently large.
\end{proof}
\begin{Cor} A countable locally finite group $G$ is finite if and only if 
\par\noindent $asdim_{AN} (G, d_G) = 0$ for every proper left invariant metric $d_G$.
\end{Cor}
\begin{proof}
It is clear that every finite group satisfies $asdim_{AN}(G, d_G) = 0$ for every metric $d_G$. \par
For the converse
firstly notice that if $(X, d_X)$ is a discrete 
metric space and $f$ is not a $0$-dimensional control function of $(X, d_X)$ then $g$ is not a $0$-dimensional control function
of $(X, d_X)$
for every function $g$ such that $g\le  f$ asymptotically that means there exists an $x_0$ such that $g(x) \le f(x)$ for every $x \ge x_0$. 
Therefore by previous theorem 
if $G$ is non finite we can take a metric $d_G$ such that the function $f(x) = x^2$ is not a $0$-dimensional control function of $(G, d_G)$.  
Then any linear function can not be a
$0$-dimensional control function of $(G, d_G)$.
\end{proof}

The remainder of this section will be focused on finding locally finite groups with non zero but still positive Assouad-Nagata dimension.
Such groups will be of the form $G = \bigoplus_{i=0}^{\infty} G_i$
with $\{G_i\}_{i\in\NN}$ some sequence of finite groups and $G_0 = \{0\}$.
Our first step consists in defining a nice proper left invariant metric in $G$ using a sequence of
 proper left invariant metrics $\{d_{G_i}\}_{i\in \NN}$ of  $\{G_i\}$. Suppose $d_{G_i}(x, y) \ge 1$ for every two different elements of $G_i$.
To build such metric we take a sequence of positive numbers
$\{s_i\}_{i\in \NN}$ such that $s_1 \ge 1$ and $s_i \ge
s_{i-1}\cdot diam(G_{i-1})+1$ and we define the map $k: G\to \NN\cup
\{0\}$ by $k(g) = \max\{i| \pi_i(g) \ne 1_{G_i}\}$. The functions
$\pi_i: G \to G_i$ are here the canonical projections. In this situation we can construct a proper norm as
following:

\begin{Lem}\label{ProperLeftInvariantMetricInG}
Let $G = \bigoplus_{i =0}^{\infty} G_i$ with $G_0 = \{0\}$ and
$\{G_i, d_{G_i}\}_{i\in \NN}$ be a sequence of finite groups with proper left invariant metrics $d_{G_i}$.
Let $\{s_i\}_{i\in \NN}$ be
the sequence of finite numbers as above.
Then the map $\|\cdot\|_G: G \to \RR_+$ defined by $\|g\|_G = s_{k(g)}\cdot \|\pi_{k(g)}(g)\|_{G_i}$ satisfies:
\begin{enumerate}
\item $\|\cdot\|_G$ is a proper norm in $G$.
\item If $g \in G_i$ then $\|j_i(g)\|_G = s_i \|g\|_{G_i}$ with $j_i: G_i \to G$ the canonical inclusion.
\end{enumerate}
\end{Lem}
\begin{proof}
(2) is obvious by definition of $\|\cdot\|_G$ so let us prove (1).
We will check firstly that $\|\cdot\|_G$ defines a norm in $G$. We
have to prove that $\|\cdot\|_G$ satisfies the three first
conditions of a norm. We are using the convention $k(g)
= 0$ if and only if $\pi_i(g) = 1_{G_i}$ for every $i$. This implies
$g = 1_G$. Hence the first
condition can be easily derived from the fact that each
$\|\cdot\|_{G_i}$ is a (proper) norm. The second condition is
trivial. For the third one let $g,h \in G$ and assume without loss
of generality that $k(g) \ge k(h)$. As $\pi_i(g\cdot h) =
\pi_i(g)$ for every $i > k(h)$  the case $k(g) >k(h)$ is obvious.
Consider now that $k(g) = k(h)$. In such case we have
$\pi_i(g\cdot h) = 1_{G_i}$ for every $i > k(g)$. There are two
possibilities:
\begin{enumerate}
\item $\pi_{k(g)}(g \cdot h) \ne 1_{G_{k(g)}}$. Then as each $\|\cdot\|_{G_i}$ is a norm we get:
\[\|g\cdot h\|_G = s_{k(g)}\cdot \|\pi_{k(g)}(g)\cdot \pi_{k(g)}(h)\|_{G_{k(g)}} \le\] \[s_{k(g)}\cdot(\|\pi_{k(g)}(g)\|_{G_{k(g)}} + \|\pi_{k(g)}(h)\|_{G_{k(g)}}) =
\|g\|_G + \|h\|_G\]
\item $\pi_{k(g)}(g\cdot h) = 1_{G_{k(g)}}$. It implies trivially $\|g\cdot h\|_G < s_{k(g)} \le \|g\|_G + \|h\|_G$.
\end{enumerate}
Finally to prove it is a proper norm  given $K >0$ we take an $s_i$ so that $K < s_i$. Hence the number of elements $g \in G$ with norm less than or equal $K$ will be
bounded by
$\Pi_{j=1}^{i-1} |G_j|$ with $|G_j|$ the cardinality of $G_j$.
\end{proof}

The proper left invariant metric $d_G$ associated to this norm will be called the {\it quasi-ultrametric generated} by $\{d_{G_i}\}_{i\in \NN}$. The reason of this name is
shown in next lemma.
\begin{Lem} Let $G = \bigoplus_{i=0}^{\infty}G_i$ be
the group defined above with $d_G$ the proper left invariant metric of the previous lemma. Then for every $g_1, g_2, g_3 \in G$
such that $k(g_i) \ne k(g_j)$ with  $i, j = 1, 2, 3$  we have:
\[d_G(g_1, g_2) \le \max\{d_G(g_2, g_3), d_G(g_1, g_3)\}\]
\end{Lem}
\begin{proof}
If $k(g_1) > k(g_2)$ and $k(g_1) >  k(g_3)$ then we will have $\pi_{k(g_1)}(g_1^{-1}\cdot g_j) = \pi_{k(g_1)}(g_1^{-1})$ with $j = 2, 3$, it will imply that
$d_G(g_1, g_j) = \|g_1^{-1}\cdot g_j\|_G = \|g_1^{-1}\|_G$. Hence $d_G(g_1, g_2) = d_G(g_1, g_j)$. Now suppose $k(g_2) > k(g_1)$ and $k(g_2) > k(g_3)$.
Applying the same reasoning we get:
$d_G(g_1, g_2) = d_G(g_2, g_3)$. Finally if $k(g_3) > k(g_1)$ and $k(g_3)> k(g_2)$ we obtain by an analogous reasoning that
$d_G(g_1, g_2) < s_{k(g_3)} \le d_G(g_3, g_2)$.
\end{proof}

We can estimate the Assouad-Nagata dimension of $G$ from the
Assouad-Nagata dimension of each $G_i$.

\begin{Lem}\label{ControlledForAllControlledForOne}
Let $G = \bigoplus_{i=0}^{\infty} G_i$ where $G_0 = \{0\}$  and
$\{G_i, d_{G_i}\}_{i \in \NN}$ is a sequence of finite groups with
$d_{G_i}$ a proper left invariant metric. Let $d_G$ be the quasi-ultrametric generated by $\{d_{G_i}\}_{i\in \NN}$.
If there is a
constant $C \ge 1$ such that for every $s\in (1, diam(G_i)]$ there exists a cover $\UU = \{\UU_0, ..., \UU_n\}$ of $(G_i, d_{G_i})$ so that the $s$-scale connected
components of each $\UU_j$ are $C\cdot s$-bounded then $asdim_{AN}(G, d_G) \le n$.
\end{Lem}
\begin{proof}
Let $s\in (s_i, s_{i+1}]$ and let $\UU = \{\UU_0, ..., \UU_n\}$ be a cover of $(G_i, d_{G_i})$
such that the $\frac{s}{s_i}$-scale connected components of each $\UU_j$ are $C\cdot \frac{s}{s_i}$ - bounded. If $s\in (s_i\cdot diam(G_i), s_i\cdot diam(G_i)+1]$ we take
as $\UU_j = G_i$ for every $j = 0, ..., n$.
Define the cover $\VV =\{\VV_0,..., \VV_n\}$ of $(G, d_G)$ by the property $g \in \VV_j$
if and only if $\pi_i(g) \in \UU_j$. Let us prove that the $s$-scale connected components of $\VV_j$ are $C\cdot s$-bounded.
Given an $s$-scale chain $x_1, x_2,...,x_m$ of $\VV_j$, define the associated chain $y_1,..., y_m$ by $y_r = x_1^{-1}\cdot x_r$. By proving
$\|y_m\|_G\le C\cdot s$ we will complete the result. \par
Firstly notice that $y_1 = 1_G$. Now suppose it is true
for some $r < m$ that $\pi_j(y_r) = 1_{G_j}$ for every $j > i$ then as a consequence of the fact
$d_G(y_r, y_{r+1}) < s \le s_{i+1}$ we will get $\pi_{j}(y_{r+1}) = 1_{G_j}$ for every $j >i$. It shows that $\pi_j(y_r) = 1_{G_j}$ for every $j>i$ and each $r = 1, ...,m$.
By construction of the metric $d_G$ we obtain $d_{G_i}(\pi_i(y_r), \pi_i(y_{r+1})) < \frac{s}{s_i}$ and then $\pi_i(y_1),...,\pi_i(y_m)$ is an $\frac{s}{s_i}$-scale chain
of $\pi_i(x_1^{-1})\cdot \UU_j$. Hence it will be $C\cdot \frac{s}{s_i}$ -bounded. Combining this and the fact $\pi_i(y_1) = 1_{G_i}$
we finally get that $\|\pi_i(y_m)\|_{G_i}\le C \cdot\frac{s}{s_i}$ what implies $\|y_m\|_G
\le C\cdot s$.
\end{proof}

So in order to get the main theorem of this section we have to find nice finite groups with proper left invariant metrics that satisfy previous lemma and
guarantee that $(G, d_G)$ has
non zero Assouad-Nagata dimension. This is the aim of next lemma. \par

\begin{Lem}\label{CisEstimatedForZk}
Let $n$ be a fixed natural number. There exists a constant $C_n\ge
1$ such that for each $s\in \RR_+$, every
finite group $(\ZZ_k^n, d_k^n)$ with $k >1$ and $d_k$ the
canonical word metric of $\ZZ_k$ has a cover $\UU = \{\UU_0,...,
\UU_n\}$ where the $s$-scale connected components of each $\UU_i$
are $C_n\cdot s$-bounded.
\end{Lem}
\begin{proof}
Define $r$ as the integer part of $\frac{k}{2}$. In this proof we will see the groups $\ZZ_k$ as:
\[\ZZ_k =
\{-r, -r+1,...-1, 0, 1,..., r-1, r\}.\] As a trivial consequence
of the results of \cite{Lang} it can be showed that $dim_{AN}
\ZZ^n \le n$ so let $D(s) = C_n'\cdot s$ be an $n$-dimensional
control function of $\ZZ^n$ with $C_n' \ge 1$. Fix $s$ a positive
number and take a cover $\VV =\{\VV_0,...,\VV_n\}$ of $\ZZ^n$ such
that the $s$-scale connected components of each $\VV_i'$ are
$C_n'\cdot s$-bounded. Define $\UU'= \{\UU_0',...,\UU_n'\}$ a
cover in $I_r^n = \{0, 1,...,r\}^n$ by the rule $\UU_i' = \VV_i
\cap I_r^n$. Notice that the restriction of $d_k^n$ to $I_r^n$
coincides with the $l_1$-metric of $I_r^n$. Let $\pi_i: \ZZ_k^n
\to \ZZ_k$ be the canonical projection over the ith coordinate.
For each subset $\lambda$ of $\{1, 2, ..., n\}$ take the
automorphism $p_{\lambda}: \ZZ_k^n \to \ZZ_k^n$ given by
$\pi_i(p_{\lambda}(x)) = \epsilon_{(\lambda, i)}\cdot \pi_i(x)$
with $\epsilon_{(\lambda, i)} = 1$ or $\epsilon_{(\lambda, i)}  =
-1$ depending on $i \in \lambda$ or $i \not\in \lambda$. Define
the cover $\UU =\{\UU_0, ...\UU_n\}$ of $\ZZ_k^n$  by the rule $x
\in \UU_i$ if and only if there exists a $\lambda$ such that
$p_{\lambda}(x) \in \UU_i'$. Let us estimate the diameter of the
$s$-scale connected components of $\UU_i$. Firstly we claim that
if $L$ and $M$ are two different $s$-scale connected components of
$\UU_i'$ then $d(p_{\lambda_1}(L), p_{\lambda_2}(M)) \ge s$ for
every $\lambda_1$ and $\lambda_2$ subsets of $\{1, 2,..., n\}$.
Suppose on the contrary that there exists $x \in L$ and $y \in M$
such that $d(p_{\lambda_1}(x), p_{\lambda_2}(y))< s$. But if $x =
(x_1, ..., x_n)$ and $y = (y_1,...,y_n)$ then:
\[d(x, y) = |x_1- y_1|+...+|x_n-y_n| \le \sum_{i=1}^n d_k(\epsilon_{(\lambda_1, i)}\cdot x_i, \epsilon_{(\lambda_2, i)} \cdot y_i) < s\]
A contradiction with the fact that $L$ and $M$ are different $s$-scale connected components. \par
From this we deduce that if $L'$ is an $s$-scale connected component of $\UU_i$ then there exists an $s$-scale connected component $L$ of $\UU_i'$ and some
subsets $\lambda_1,...,\lambda_m$ of $\{1, ..., n\}$ so that $L' = \bigcup_{j=1}^m p_{\lambda_j}(L)$. Each $p_{\lambda_j}(L)$ is $s$-scale connected and $C_n'\cdot s$-bounded
as the maps $p_{\lambda}$ are isometries.  It is clear that $m \le 2^n$. Hence we get:
\[diam(L') = diam( \bigcup_{j=1}^m p_{\lambda_j}(L))\le \sum_{i=1}^m(C_n'\cdot s + s) \le 2^n\cdot(C_n'+1)\cdot s\]
Therefore the cover $\UU =\{\UU_0,...,\UU_n\}$ satisfies the
requirements of the lemma with $C_n = 2^n\cdot(C_n'+1)$.
\end{proof}
\begin{Thm}\label{MainTheorem} For each $n \in \NN \cup \{\infty\}$ there exists a locally finite group
$G^n$ with a proper left invariant metric $d_n$ such that
$asdim_{AN} (G^n,d_n) =n$.
\end{Thm}
\begin{proof}
Case $n < \infty$. Suppose $n$ fixed. Let $\{k_i\}_{i\in \NN}$ be
an increasing sequence of natural numbers with $k_1 >1$ and let
$\{r_i\}_{i\in \NN}$ be the sequence given by the integer part of
$\frac{k_i}{2}$. Take the finite groups $(G_i^n = \ZZ_{k_i}^n,
d_{G_i^n})$ with $d_{G_i^n}$ the canonical word metric and let
$(G, d_G)$ be the group $G = \bigoplus_{i = 1}^{\infty} G_i^n$
with $d_G$ the quasi-ultrametric generated by $\{d_{G_i^n}\}$. It
is clear that if we embed the subsets $\{0, 1,...,r_i\}^n$ of each
$G_i^n$ in $G$ we get a sequence of $n$-dimensional dilated
cubes of increasing size so applying
\ref{SequenceImpliesNagata} we obtain $asdim_{AN}(G, d_G) \ge n$.
On the other hand by lemmas \ref{ControlledForAllControlledForOne}
and \ref{CisEstimatedForZk} we conclude $asdim_{AN}(G, d_G) \le
n$.\par Case $n = \infty$. Just take the group $G = \bigoplus_{i=
1}^{\infty} G_i^i$ with $G_i^i$ the groups defined in the previous
case and construct the quasi-ultrametric $d_G$ generated by
$\{d_{G_i^i}\}_{i\in \NN}$. Applying an analogous reasoning as
above we get that $asdim_{AN}(G, d_G) \ge n$ for each $n$.
\end{proof}
\begin{Rem}
It was asked in \cite{Brod-Dydak-Levin-Mitra} about the possibility of defining the asymptotic
Assouad-Nagata dimension of arbitrary groups $G$ as the supremum of the asymptotic Assouad Nagata dimensions of its finitely generated subgroups $H$.
From the previous results
we can deduce that this approach does not work well. Even if we assume that the asymptotic Assouad-Nagata dimension of $G$ is finite.
\end{Rem}
\begin{Cor}\label{GeneralizationToCountable} For every $n, k$ with $n \in \NN$ and $k\in \NN \cup \{\infty\}$ there exists a countable abelian group $(G^{(n, k)}, d_{(n, k)})$ with $d_{(n,k)}$
a proper left invariant metric such that
$asdim (G^{(n, k)}, d_{(n, k)}) = n$ but $asdim_{AN}((G^{(n, k)}, d_{(n, k)}) = n+k$.
\end{Cor}
\begin{proof}
Fix $n$ and $k$ as in the hypothesis and take the group $(G^k, d_k)$ as in theorem \ref{MainTheorem}. Define the group $G^{(n, k)} = \ZZ^n\oplus G^k$ with
the proper left invariant metric $d_{(n, k)}$ given by $d_{(n,k)}((x_1, y_1), (x_2, y_2)) = \|x_1 -x_2\|_1 + d_k(y_1, y_2)$. Multiplying an $n$-dimensional
dilated cube of $\ZZ^n$ with a $k$-dimensional dilated cube of $G^k$ we will get an $n+k$-dimensional dilated cube in $G^{(n, k)}$.
Applying \ref{SequenceImpliesNagata}
we deduce $\asdim_{AN}(G^{(n, k)}, d_{(n,k)}) \ge n + k$. The other inequalities follow easily from the 
by the subadditivity of the asymptotic dimension and the Assouad-Nagata dimension with
respect to the cartesian product(see for example \cite{Brod-Dydak-Levin-Mitra}) wand the well known fact $\asdim(\ZZ^n, d_1) = n$. 
\end{proof}

\begin{Problem} Does any countable group $G$ of finite asymptotic dimension satisfy the following condition: There exists a proper left invariant metric $d_G$ such that 
$asdim(G, d_G) < asdim_{AN}(G, d_G) < \infty$?
\end{Problem}

\end{document}